\documentclass[11pt]{article}         %[neutral.tex 27.08.07]
\usepackage{amsfonts,amssymb}
\makeatletter
\def\rdraft{\pagestyle{myheadings}
            \topmargin=30pt\headheight=10pt\headsep=3pt\footskip=20pt
            \textheight=10.5truein \textwidth=7.5truein
            \parindent=8pt \voffset=-1truein
            \ifcase \@ptsize \hoffset=-1.5truein \or \hoffset=-1.35truein
                    \or \hoffset=-1.15truein \fi}
\def\quality{\textheight=240mm \textwidth=160mm \topmargin=0Truein
             \ifcase \@ptsize \hoffset=-23mm
                     \or \hoffset=-20mm \or \hoffset=-15mm \fi}
\makeatother %\rdraft
\quality

\def\bline(#1,#2)(#3,#4)(#5){\put(#1,#2){\line(#3,#4){#5}}}  %straight line

\newcommand\mlbscale{1pt} %to change: \renewcommand\mlbscale{1.3pt}
\newif\iffigs\figstrue %\newif\iffigs\figsfalse -- to fake figures

\def\bfig(#1,#2)#3#4{\begin{figure} \begin{center}
    \framebox{\setlength{\unitlength}{\mlbscale}
       \iffigs \begin{picture}(#1,#2) #3 \end{picture}
       \else \begin{picture}(60,10)(0,0)
                   \put(0,0){\framebox(60,10){Figure}} \end{picture} \fi}
    \end{center} \caption{#4} \end{figure}}
\def\Bfig(#1,#2)#3#4{\begin{figure} \begin{center}
    \setlength{\unitlength}{\mlbscale}
       \iffigs \begin{picture}(#1,#2) #3 \end{picture}
       \else \begin{picture}(60,10)(0,0)
                   \put(0,0){\framebox(60,10){Figure}} \end{picture} \fi
    \end{center} \caption{#4} \end{figure}}
\def\bpic(#1,#2)#3{\setlength{\unitlength}{\mlbscale}
    \begin{picture}(#1,#2) #3 \end{picture}}
\def\n{\noindent}    \def\ep{\varepsilon}

\def\function#1{\left\{\!\!\!\begin{array}{ll} #1 \end{array} \right.}
\def\bea#1{\begin{eqnarray*} #1 \end{eqnarray*}} \def\a{\!\!\!&\!\!\!\!&}
\def\beq#1#2{\begin{equation} \label{#1} #2 \end{equation}}
\def\beaq#1#2{\label{#1} \begin{eqnarray} #2 \end{eqnarray}}
\def\thname{Theorem}     \def\lmname{Lemma}      \def\prname{Proposition}
\def\dfname{Definition}  \def\crname{Corollary}  \def\rmname{Remark}

\newtheorem{theorem}{\thname}%[section]   %Numbering: Theorem--Other section
\newtheorem{lemma}{\lmname}%[section]     %{lemma}[theorem]{Lemma}   section
\newtheorem{corollary}[lemma]{\crname}   %lemma

\newtheorem{dftn}{\dfname}[section]
\newtheorem{rmrk}[lemma]{\rmname}
\newenvironment{remark}{\begin{rmrk}\rm}{\end{rmrk}}     %lemma

            \catcode`@=11 \@addtoreset{equation}{section} \catcode`@=12

\def\proof{\smallskip \noindent {\bf Proof. \ }}
\newcommand\emptysquare{\ \ $\Box$}
\def\qed{\hfill\emptysquare\linebreak\smallskip\par}
\def\*#1{#1^*}    \def\map{T}  
  
\def\noprint#1{}    
     \def\cM{\cal M}

   \def\t#1{\tilde#1}

\def\toas#1{\stackrel{#1}{\longrightarrow}}
\def\blim#1#2{\if #1+ \limsup_{#2} \else {\if #1- \liminf_{#2} \else
              \lim_{#2}\left(\begin{array}{l}\sup\\
              \inf\end{array}\right) \fi} \fi} %\blim\pm{t\to0}

  \def\la{\lambda}  \def\*#1{{#1^*}}
   
\def\lal{\la_\alpha}  \def\jal{{j_\alpha}}

\def\?#1{} % comments

\begin{document}%%%-------------------------------------------------%%%%

\title{Finite rank approximations of expanding \\
       maps with neutral singularities}
\author{Michael Blank\thanks{
        Russian Academy of Sci., Inst. for
        Information Transm. Problems,
        and Laboratoire Cassiopee UMR6202, CNRS, ~
        e-mail: blank@iitp.ru}
        \thanks{This research has been partially supported
                by Russian Foundation for Fundamental Research, CRDF
                and French Ministry of Education grants.}
       }
\date{August 27, 2007} %\today} %April 25, 2007}
\maketitle

%\vskip-2cm

\begin{abstract} For a class of expanding maps with neutral
singularities we prove the validity of a finite rank
approximation scheme for the analysis of Sinai-Ruelle-Bowen
measures. Earlier results of this sort were known only
in the case of hyperbolic systems.
\end{abstract}

\n{\bf AMS Subject Classification}: Primary 37M25, 37A40;
Secondary 37A30, 37A50, 37C30.

\n{\bf Key words}: dynamical system, SRB measure, transfer operator,
Ulam's approximation.

%\bigskip%
%\n{\bf Key words}: dynamical system, SRB measure.

%\bigskip%
%

\section{Introduction}\label{s:intro}

%{\em Finite rank approximations of chaotic dynamical systems with neutral singularities}

In 1960 S. Ulam \cite{Ul} has formulated a hypothesis about the
possibility of an approximation of an action of a chaotic
dynamical system by means of a sequence of finite state Markov
chains. He even proposed the simplest scheme for such an
approximation which can be described in modern terms as follows.
Let $\map$ be a map from a Lebesgue compact space $(X,m)$ equipped
with a metric $\rho$ into itself. Iterations $\map,
\map^2\!\!:=\!\!\map\circ\map, \map^3, \dots$ of the map $\map$ define
a discrete time {\em dynamical system} on $X$. One extends the
action of the map $\map$ to the set of probabilistic measures
(generalized functions) on $X$ according to the formula: %
$$ \*\map\mu(A) := \mu(\map^{-1}A) $$
for any Borel set $A\subseteq X$. We shall refer to $\*\map$
as a {\em transfer-operator} corresponding to the dynamical system
$(\map,X)$. Let $\Delta:=\{\Delta_i\}$ be a finite measurable
partition of $X$ with the diameter $\delta$. Consider an operator
acting on probabilistic measures (generalized functions):
$$ \*{Q_\Delta}\mu(A)
:= \sum_i\frac{m(A\cap\Delta_i)}{m(\Delta_i)}~\mu(\Delta_i) .$$

In this terms the Ulam's approximation can be written as a
superposition of the operators $\*{Q_\Delta}$ and $\*\map$,
and his hypothesis says that for a ``good'' enough map and
``good'' enough partitions $\Delta$ statistical properties of
the original dynamical system can be obtained from the limit
properties of the ``spatially discretized'' transfer operators
$\*{\map_\Delta}:=\*{Q_\Delta}\*\map$ when the partition diameter
vanishes. In particular, the so called Sinai-Ruelle-Bowen (SRB)
measure of the dynamical system corresponds to the limit of
leading eigenfunctions of the operator $\*{\map_\Delta}$
considered as a linear operator in a suitable Banach space of
signed measures (generalized functions). %
Recall that the {\em SRB measure} is a probabilistic measure
$\mu_\map$ satisfying the property that there is an open subset
$A\subseteq X$ such that %
$\frac1n\sum_{k=0}^{n-1}{\*\map}^k\mu\toas{n\to\infty}\mu_\map$ %
for any probabilistic measure $\mu$ absolutely continuous with
respect to the reference measure $m$ and having the support on
$A$. This version of the SRB measure is often called a natural or
physical measure. We refer the reader to \cite{Bl-mon} for
detailed discussions of SRB measures and their properties.

From a numerical point of view the operator $\*{\map_\Delta}$
is equivalent to a transition matrix $P=(p_{ij})$ of a finite
state Markov chain with transition probabilities
$p_{ij} = m(\Delta_i\cap\map^{-1}\Delta_j)/m(\Delta_i)$.
Therefore its complete analysis on a computer is a routine
procedure (see \cite{dellnitz} for details).

The main problem with the analysis of the Ulam type approximation
is how to connect the dynamics of Markov chains defined by the
approximation with the original dynamics. One is tempted here to
adapt the partition $\Delta$ to geometric properties of the map,
in particular, to use the so called Markov partitions (see e.g.
\cite{Fr}). This idea simplifies the analysis a lot making it
similar to classical symbolic dynamics. Unfortunately in practice
the usefulness of the adapted partitions is limited by the
observation that usually such partitions can be found only
numerically. Therefore small errors are inevitable and they may
lead to even worse accuracy compared to a generic partition
(see \cite{BK97} for details).

A natural next step here is to analyze connections between the
complete spectrum of the original transfer-operator and the limit
of the spectra of the perturbed transfer-operators. It turns out
that for a broad class of dynamical systems having some
hyperbolicity properties (piecewise expanding maps \cite{Li, BK95,
BK97, Bl-mon}, Anosov torus diffeomorphisms \cite{BKL}, random
maps \cite{Bl3}) one might show that both the corresponding
transfer-operator and its perturbation are quasi-compact (i.e. is
a sum of a compact operator and a finite dimensional projector).
Using this property it is possible to prove that the part of the
spectra corresponding to isolated eigenvalues indeed, satisfies
the above mentioned hypothesis (see also \cite{Bl-mon,dellnitz,
Fr} for the discussion of numerical realizations of finite rank
approximations).

Strictly speaking even for a very ``good'' hyperbolic dynamical
system some additional assumptions are necessary to prove the
hypothesis for all isolated eigenvalues. Surprisingly, a similar
statement about the leading eigenfunction turns out to be
extremely robust. In fact, the only known counterexample (see
below) is not only discontinuous but this discontinuity occurs at
a periodic turning point (compare to instability results about
general random perturbations in \cite{BK97}).

\begin{lemma}\label{l:counter-or-ulam} \cite{Bl3} The map
$$ \map x := \function{
            \frac{x}4+\frac12 &\mbox{if }\; 0\le x <\frac5{12}\\
            -2x +1            &\mbox{if }\; \frac5{12}\le x < \frac12\\
            \frac{x}2+\frac14 &\mbox{otherwise} .} $$
from the unit interval into itself is uniquely ergodic, but the
leading eigenvector of the Ulam approximation
$\*{Q_{\Delta}}\*\map$ corresponding to the partition into
intervals of the same length does not converge weakly to the only
$\map$-invariant measure.
\end{lemma}

Up to now there were no mathematical results corresponding to the
situation when the transfer-operator has no isolated eigenvalues.
Despite the conventional techniques mentioned above no longer
works in this case we shall prove the stability of the leading
eigenfunction (which corresponds to the SRB measure) for some
nonhyperbolic systems.

Consider a family of expanding maps with neutral singularities.
A typical example of this type is the so called Manneville-Pomeau
map $\map_\alpha x:= x+x^{1+\alpha} ({\rm mod}~ 1)$
from the unit interval $X:=[0,1]$ into itself with $\alpha>0$. The
interest to such systems is to a large extent due to the fact that
they model the so called intermittency phenomenon \cite{PM}.
It is well known (see e.g. \cite{Th, LVS}) that map
$\map_\alpha$ possesses the only one SRB measure $\mu_\alpha$
and that this measure is absolutely continuous (but has an
unbounded density) with respect to the Lebesgue measure $m$
if $0<\alpha<1$, while for $\alpha>1$ it
coincides with the Dirac measure at the origin $\*{1_{\{0\}}}$.

Let $\Delta$ be a partition with diameter $\delta$ into
(unnecessary equal) intervals $\{\Delta_i\}$ satisfying the
property $m(\Delta_i)/m(\Delta_j)\le{K}<\infty~~\forall i,j$ and
let $\Delta_1\in\Delta$ be the interval containing the origin. The
following result demonstrates that the Ulam scheme of finite rank
approximations of this nonhyperbolic map is correct.

Recall that a Markov chain is uniquely ergodic if it has a unique
stationary probability distribution.

\begin{theorem}\label{t:main} For any $\alpha\ge0$ and small
enough $0<\delta\ll1$ the Markov chain generated by the
transfer operator $\*{Q_\Delta} \*{\map_\alpha}$ is uniquely
ergodic and its unique invariant distribution $\mu_\Delta$
satisfies the relations: %
\begin{enumerate}
\item[(a)] $\mu_\Delta(\Delta_1)\le C\delta^{1-\alpha}
            \quad\forall\alpha\ge0$,

\item[(b)] $\mu_\Delta(\Delta_1)/\delta\toas{\delta\to0}\infty
            \quad\forall\alpha>0$,

\item[(c)] $\mu_\Delta\toas{\delta\to0}\*{1_{\{0\}}}
            \quad\forall\alpha>1$.

\end{enumerate}
\end{theorem}

Due to the nonhyperbolicity of the map $\map_\alpha$ the operator
approach discussed above no longer works here while methods used
in the analysis of maps with neutral singularities do not work
with highly discontinuous densities unavoidable due to the
action of the projection operator $\*{Q_\Delta}$. Therefore we
develop a completely new approach based on the analysis of the
action of the corresponding transfer operators on ``monotonic
measures'' $\mu$ defined by the property that
$$ \mu(A) \ge \mu((A+x)\cap X) $$
for any interval $A\subset X:=[0,1]$ and any number $x\in X$,
where $A+x:=\{a+x:~~ a\in A\}$.

In fact, these results hold for a much more general class of
expanding maps with neutral singularities and we shall discuss
sufficient conditions for them in Section~\ref{s:generalization}.

\section{Action of transfer operators on monotonic measures}
\label{s:monotonic}

Denote by $\cM$ the set of all monotonic probabilistic
measures on $X$.

\begin{lemma}\label{l:monotonic}
Each element $\mu\in{\cM}$ can be uniquely represented as a
weighted sum of the Dirac measure at zero and an absolutely
continuous measure (with respect to $m$) having a monotonous
non-increasing density.
\end{lemma}
\proof Recall that a measure $\mu$ is absolutely continuous if
and only if for any $\ep>0$ there exists $\delta>0$ such that for
any finite collection $A$ of nonintersecting intervals of total
length $m(A)\le\delta$ one has $\mu(A)\le\ep$. Let $\mu\in{\cM}$.
Assume that for some $b\in(0,1)$ and $0<\delta<1-b$ the interval
$B:=[b,b+\delta]$ satisfies the inequality $\mu(B)\le\ep$. Then
intervals from any finite collection $A\subset[b,1]$ of
nonintersecting intervals of total length $m(A)\le\delta$ may be
shifted to the left (preserving their respective lengths) such
that the shifted collection $A'$ will be still a collection of
nonintersecting intervals but belonging to the interval $B$. This
is always possible since $m(A)\le\delta=m(B)$. The monotonicity
property of the measure $\mu$ implies $\mu(A')\ge\mu(A)$. On the
other hand, $A'\subset B$ and thus $\mu(A')\le\mu(B)\le\ep$. Thus
the restriction of $\mu$ to the interval $[b,1]$ is an absolutely
continuous measure.

Let us show that for any $b\in(0,1)$ and any $\ep>0$ there
exists $0<\delta<1-b$ such that $\mu(B)\le\ep$. Assume that this
is not the case and there exists a pair $b,\ep$ such that for any
$0<\delta<1-b$ we have $\mu([b,b+\delta])>\ep$. Then by definition
$$ \mu([b-n\delta,b+\delta-n\delta])\ge \mu([b,b+\delta]) \ge \ep $$
for any positive integer $n\le N:= \lfloor b/\delta\rfloor$, where
$\lfloor \cdot \rfloor$ stands for the integer part. Thus %
\bea{ 1 \a\ge \mu([0,b+\delta])
       \ge \sum_{n=0}^N \mu([b-n\delta,b+\delta-n\delta]) \\
     \a\ge (N+1) \mu([b,b+\delta]) \ge (b/\delta-1)\ep
     \toas{\delta\to0}\infty .}%
We came to the contradiction.

Therefore the only place where a singular component of a
monotonic measure $\mu$ may appear is the origin. On the other
hand, the Dirac measure at the origin clearly satisfies the
monotonic property which proves the representation in the form of
the weighted sum. The uniqueness of this representation follows
from a general result about the decomposition of a measure into
singular and absolutely continuous components.

It remains to prove that the density $f:=d(\mu-\mu(\{0\}))/dm$ of
the absolutely continuous component is a monotonous
non-increasing function. By the monotonicity of the measure $\mu$
for any interval $A\in(0,1]$ and any $x>0$ such that
$y+x\le1~~\forall y\in A$ we have %
\bea{ \mu(A) \a= \int_A f(y) dm(y)
             \ge \mu(A+x) \\
             \a= \int_{A+x} f(y) dm(y)
             \ge \int_A f(y+x) dm(y) .} %
Since $A,x$ are arbitrary this implies that the density $f$ is a
non-increasing function on a subset $Y\subseteq X$ of full
Lebesgue measure. Redefining $f$ on the complement to this set as
$f(x):=\inf\limits_{Y\ni y<x}f(y)$ we obtain a representative of
the same $L^1(m)$-equivalence class for which the monotonicity
holds everywhere. \qed

\begin{corollary}\footnote{The author is grateful to an anonymous
    referee  for this characterization of monotonic measures.} %
A measure $\mu\in\cM$ iff $\mu([0,x])$ is a convex function  on $x\in X$.
\end{corollary}

The definition of the monotonic measure makes it possible to
compare its values directly only on intervals of the same length.
The following result extends this property for intervals of
different lengths and technically is one of the key ingredients
of our approach.

\begin{lemma}\label{l:key} For any two nonempty intervals
$A,B\subseteq X$ such that %
\beq{e:int-ineq}{\inf\{a\in A\}\le\inf\{b\in B\}, \qquad
                 \sup\{a\in A\}\le\sup\{b\in B\}} %
and any monotonic measure $\mu$ we have %
\beq{e:den-monotone}{\frac{\mu(A)}{m(A)} \ge \frac{\mu(B)}{m(B)} .}
\end{lemma}

\proof To simplify notation for a Borel set $A\subseteq X$ we
denote $|A|:=m(A)$, $\mu_A:=\mu(A)/|A|$ if $|A|>0$ and $\mu_A=0$
otherwise. By $f$ we denote the density of the absolutely
continuous component of the measure $\mu$. We say also that $A\le
B$ if the nonempty intervals $A,B\subseteq X$ satisfy
(\ref{e:int-ineq}).

Using this notation our claim can be written as $\mu_A\ge\mu_B$
whenever $A\le B$ and means basically that the average density
decays when the interval of averaging moves to the right.

By the additivity of measures we get
$$ \mu_A
 = \frac1{|A|} \left( |A\setminus B|\cdot\mu_{A\setminus B}
                    + |A\cap B|\cdot\mu_{A\cap B} \right) ,$$
$$ \mu_B
 = \frac1{|B|} \left( |B\setminus A|\cdot\mu_{B\setminus A}
                    + |A\cap B|\cdot\mu_{A\cap B} \right) .$$
Thus %
\bea{\a |A|\cdot|B|\cdot(\mu_A - \mu_B) \\
 \a= |B|\cdot|A\setminus B|\cdot\mu_{A\setminus B}
 + |B|\cdot|A\cap B|\cdot\mu_{A\cap B} \\
 \a\quad - |A|\cdot|B\setminus A|\cdot\mu_{B\setminus A}
 - |A|\cdot|A\cap B|\cdot\mu_{A\cap B} .}%
Observe that %
\bea{ \a |B|\cdot|A\setminus B|\cdot\mu_{A\setminus B}
 - |A|\cdot|A\cap B|\cdot\mu_{A\cap B} \\
 \a= |B\setminus A|\cdot|A\setminus B|\cdot\mu_{A\setminus B}
 + |A\cap B|\cdot|A\setminus B|\cdot\mu_{A\setminus B} \\
 \a\quad - |A\setminus B|\cdot|A\cap B|\cdot\mu_{A\cap B}
 - |A\cap B|^2\mu_{A\cap B} \\ %
 \a= |B\setminus A|\cdot|A\setminus B|\cdot\mu_{A\setminus B}
 + |A\setminus B|\cdot|A\cap B|
                 \cdot(\mu_{A\setminus B} - \mu_{A\cap B}) \\
 \a\quad - |A\cap B|^2\mu_{A\cap B} .}%
Similarly %
\bea{ \a |A|\cdot|B\setminus A|\cdot\mu_{B\setminus A}
 - |B|\cdot|A\cap B|\cdot\mu_{A\cap B} \\
 \a= |A\setminus B|\cdot|B\setminus A|\cdot\mu_{B\setminus A}
 + |B\setminus A|\cdot|A\cap B|
                 \cdot(\mu_{B\setminus A} - \mu_{A\cap B}) \\
 \a\quad - |A\cap B|^2\mu_{A\cap B} .}%
Therefore %
\bea{\a |A|\cdot|B|\cdot(\mu_A - \mu_B) \\
 \a= |B\setminus A|\cdot|A\setminus B|
              \cdot(\mu_{A\setminus B} - \mu_{B\setminus A})
 + |A\setminus B|\cdot|A\cap B|
              \cdot(\mu_{A\setminus B} - \mu_{A\cap B}) \\
 \a\quad + |B\setminus A|\cdot|A\cap B|
              \cdot(\mu_{A\cap B} - \mu_{B\setminus A}) .}%

If $|A\setminus B|\cdot|A\cap B|>0$ by the monotonicity of the
density $f$ and that $0\notin B$ (since $|A\setminus B|>0$) we have
$$ \mu_{A\setminus B} \ge \inf_{A\setminus B} f
                      \ge \sup_{A\cap B} f \ge \mu_{A\cap B}
                      \ge \inf_{A\cap B} f
                      \ge \sup_{B\setminus A} f
                      \ge \mu_{B\setminus A} $$
and hence all summands in the above sum are nonnegative,
which implies that $\mu_A\ge\mu_B$.

If $|A\setminus B|=0$ and $|A\cap B|>0$ we have $A=A\cap B$
(at least up to the endpoints) and
$$ |A|\cdot|B|\cdot(\mu_A - \mu_B)
 = |A|\cdot|B\setminus A|
              \cdot(\mu_{A} - \mu_{B\setminus A}) .$$
In this case $0\notin B\setminus A$ and making use of
the monotonicity of the density we get:
$$ \mu_A \ge \inf_A f \ge \sup_{B\setminus A} f
         \ge \mu_{B\setminus A} ,$$
which proves the claim in this case.

It remains to consider the simplest case
$|A\setminus B|>0$ and $|A\cap B|=0$. Here again
$0\notin B$ and following the same argument as above we get
$$ \mu_A \ge \inf_A f \ge \sup_B f \ge \mu_B .$$
\qed

\begin{remark} The condition (\ref{e:int-ineq}) is crucial
here and cannot be relaxed: if any of the inequalities
(\ref{e:int-ineq}) is violated, the inequality
(\ref{e:den-monotone}) might not hold as well.
Indeed, for $f(x):=2-2x, ~A=[0,1], ~ B=[1/4,1/2]$
we have $\mu_A=1/2<5/8=\mu_B$ despite
$\inf\{a\in A\}=0<1/4=\inf\{b\in B\}$.
\end{remark}

\begin{lemma}\label{l:Ulam-inv} $\*{Q_\Delta}\cM\subset\cM$.
\end{lemma}
\proof
We assume always that the intervals $\Delta_i$ are enumerated
in a natural way according to their positions, namely that
$\Delta_i\le \Delta_j$ if $i\le j$.

Observe that by the definition of the transfer operator
$\*{Q_\Delta}$ the measure $\*{Q_\Delta}\mu$ is absolutely
continuous irrespective of the measure $\mu$. Therefore for
$\mu\in\cM$ the measure $\*{Q_\Delta}\mu$ always has a density
which we denote by $\t{f}$.

Since $\t{f}_{|\Delta_i}(x)$ is a constant for all $x\in\Delta_i$
we can drop the dependence on $x$. For any $i<j$ by
Lemma~\ref{l:key} we have
$$ \t{f}_{|\Delta_i}
 = \frac1{|\Delta_i|} \*{Q_\Delta}\mu(\Delta_i)
 = \frac1{|\Delta_i|} \mu(\Delta_i)
 = \mu_{\Delta_i} \ge \mu_{\Delta_j}
 = \frac1{|\Delta_j|} \*{Q_\Delta}\mu(\Delta_j)
 = \t{f}_{|\Delta_j} .$$

This relation proves that the density of the measure
$\*{Q_\Delta}\mu$ is non-increasing, which immediately implies
that this measure is monotonic. \qed

Consider now a family of piecewise convex maps ${\cal T}$ from
the unit interval $X$ into itself such that for each map
$\map\in{\cal T}$ there is a partition $\{X_i\}$ of $X$ into
intervals (called {\em special} partition) satisfying the
following properties:
\begin{itemize}
\item $\map_{|X_i}:X_i\to\map X_i$ is a convex one-to-one map
      for each $i$.
\item $0\in \map X_i$ for each $i$.
\end{itemize}

Observe that these two assumptions imply that $\map_{|X_i}$ is
monotonous increasing. We assume also that the intervals
belonging to the special partition are ordered in a natural way,
i.e. $X_i\le X_j$ (in the sense of (\ref{e:int-ineq})) if $i<j$.
Therefore $0\in X_1$. A typical example of a map $\map\in{\cal
T}$ is represented on Fig.~\ref{f:convex-map}.

%%%%%%%%%%%%%%%%%%%%%%%%%%%%%%%%%%%%%%%%%%%%
%% Example of a piecewise convex map
\Bfig(150,150)
      {\footnotesize{
       \bline(0,0)(1,0)(150)   \bline(0,0)(0,1)(150)
       \bline(0,150)(1,0)(150) \bline(150,0)(0,1)(150)
       %\bezier{100}(0,75)(75,75)(150,75)
       %\bezier{100}(75,0)(75,75)(75,150)
       %\bezier{70}(60,0)(60,45)(60,90)
       %\put(0,75){\vector(4,1){60}}  \put(60,30){\vector(1,-2){15}}
       %\bline(75,75)(2,1)(75)
       \bezier{100}(0,0)(75,75)(150,150)
       \bezier{100}(90,0)(90,75)(90,150)
       \bezier{100}(125,0)(125,75)(125,150)
       \thicklines
       \bezier{200}(0,0)(60,20)(90,120)
       \bezier{200}(90,0)(110,20)(125,140)
       \bezier{200}(125,0)(140,0)(149,90)
       \put(-5,-8){$0$} \put(147,-8){$1$}
       \put(40,-8){$X_1$} \put(100,-8){$X_2$}
       \put(130,-8){$X_3$} \put(-7,145){$1$}
      }}
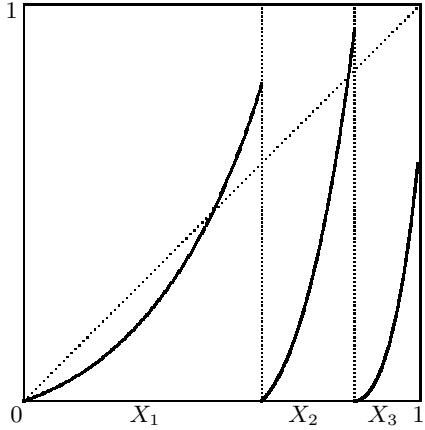
{An example of a map from the family ${\cal T}$.
\label{f:convex-map}}
%%%%%%%%%%%%%%%%%%%%%%%%%%%%%%%%%%%%%%%%%%%%%%%%%%%%%%%%%%%%%%%%%

\begin{lemma}\label{l:map-inv}
$\*\map\cM\subset\cM$ for any $\map\in {\cal T}$.
\end{lemma}
\proof We need to show that if $\mu\in\cM$ then for any pair of
nonempty intervals $A,B\subset X$ with $A\le B$ (in the sense of
(\ref{e:int-ineq})) and $|A|=|B|$ one has
$$ \*\map\mu(A) \ge \*\map\mu(B) .$$
Denoting $\map_i:=\map_{|X_i}$ we get
$$ \*\map\mu(A) = \mu(\map^{-1}A)
                = \mu(\cup_i \map_i^{-1}A)
                = \sum_i \mu(\map_i^{-1}A) $$
since $\map_i^{-1}A\cap\map_j^{-1}A=\emptyset$
for any $i\ne j$.

On the other hand, $\map_i^{-1}$ (being an inverse map to a
convex one) is a concave map and it preserves the origin for
each $i$. Fix some $i$ and consider a concave origin
preserving map $G:\map X_i\to X$.

Recall that a function is concave if for any pair of points the
straight line connecting their values lies below the graph of the
function. Since $G0=0$ the concave map $G$ is a monotone
increasing one-to-one continuous map. Thus for any interval
$A=[a,a']$ its image $GA=[Ga,Ga']$ and $m(GA)=Ga'-Ga$.

Our aim now is to show that $\mu(GA) \ge \mu(GB)$ whenever
$A\le B$, $|A|=|B|$. Since $GA\le GB$ (by the monotonicity
of $G$) for $\mu\in\cM$ by Lemma~\ref{l:key}:
$$ \frac{\mu(GA)}{m(GA)} \ge \frac{\mu(GB)}{m(GB)} .$$
Hence
$$ \mu(GA) \ge \frac{m(GA)}{m(GB)}~\mu(GB) $$
and it remains to prove only that $m(GA)\ge m(GB)$.

There are two possibilities: either $|A\cap B|=0$ or
$|A\cap B|>0$. We start with the first case, i.e.
$a<a'\le b<b'$, where $B=[b,b']$. Since $G$ is concave
and continuous the slopes of the straight lines
connecting consecutively the points $(a,Ga), (a',Ga'),
(b,Gb), (b',Gb')$ do not increase from one interval
to another, i.e. %
\beq{e:slopes}{
  \frac{|GA|}{|A|} \ge \frac{m(G([a',b]))}{m([a',b])}
                   \ge  \frac{|GB|}{|B|} }%
provided $a'<b$ (otherwise the middle term should be dropped).
Thus $|GA|\ge|GB|$. Note the similarity between (\ref{e:slopes})
and (\ref{e:den-monotone}).

Consider the second case $|A\cap B|>0$, i.e. $a<b<a'<b'$.
We have
$$ GA = [Ga,Gb) \cup [Gb,Ga'], \qquad
   GB = [Gb,Ga'] \cup (Ga',Gb'] $$
and these unions a disjoint. Therefore
$$ |GA| - |GB| = m([Ga,Gb)) - m((Ga',Gb'])) \ge 0 $$
since (according to the decrease of the slopes)
$$ \frac{m([Ga,Gb])}{m([a,b])}
\ge \frac{m([Ga',Gb'])}{m([a',b'])} $$
and
$$ m([a,b]) = |A| - m([b,a']) = |B| - m([b,a'])
           = m([a',b']) .$$
This finishes the proof that $\mu(GA)\ge\mu(GB)$ for any concave
origin preserving map $G$. Returning to the original notation we
get
$$ \*\map\mu(A) = \sum_i\mu(\map_i^{-1}A)
              \ge \sum_i\mu(\map_i^{-1}B) = \*\map\mu(B) .$$
\qed

\section{Proof of Theorem~\ref{t:main}} \label{s:proof}

It is straightforward to check that under the repeated
applications of the map $\map_\alpha$ any interval covers the
entire phase space $X$ in a finite number of iterations.
Therefore for a given finite partition $\Delta=\{\Delta_i\}$ into
intervals there exists a positive integer $n_\Delta<\infty$ such
that $\map_\alpha^{n_\Delta}\Delta_i=X$ for each $i$ and hence
the transition matrix $P=(p_{i,j})$ corresponding to the operator
$\*{Q_\Delta}\*{\map_\alpha}$ in power $n_\Delta$ is strictly
positive (i.e. all its entries are positive). Recall that
$$ p_{i,j}:=\frac{m(\map_{\alpha}^{-1}\Delta_{j}\cap\Delta_{i})}
                 {m(\Delta_{i})} .$$
Then by the Perron-Frobenius Theorem the Markov matrix $P$ has
the only one normalized left eigenvector $\pi=(\pi_i)$ with the
unit eigenvalue,  i.e.
$$ \sum_i\pi_i=1, \qquad \sum_i p_{i,k}\pi_i=\pi_k ~~\forall k .$$
Moreover, applying iteratively the matrix $P$ to any normalized
vector with nonnegative entries one converges to $\pi$. Denote a
probabilistic measure $\mu_\Delta$ by the relation
$$ \mu_\Delta(A) = \sum_i \pi_i ~m(A\cap \Delta_i) $$
for any Borel set $A\subseteq X$. Clearly this implies
$\mu_\Delta(\Delta_i)=\pi_i$ for all $i$. The uniqueness of $\pi$
and the convergence to it for any nonnegative initial vector
immediately implies that $\mu_\Delta$ is a SRB measure for the
operator $\*{Q_\Delta} \*{\map_\alpha}$.

Observe now that for any $\alpha\ge0$ the map $\map_\alpha$
(presented on Fig~\ref{f:M-P}) belongs to the family of piecewise
convex maps ${\cal T}$ defined in Section~\ref{s:monotonic}.
Therefore Lemmas \ref{l:Ulam-inv} and \ref{l:map-inv} imply that
$\*{Q_\Delta} \*{\map_\alpha}\cM\subset\cM$. On the other hand, as
we just demonstrated the measure $\mu_\Delta$ is the only
invariant measure of this process, therefore $\mu_\Delta\in\cM$.

%%%%%%%%%%%%%%%%%%%%%%%%%%%%%%%%%%%%%%%%%%%%
%% Example of a piecewise convex map
\Bfig(150,150)
      {\footnotesize{
       \bline(0,0)(1,0)(150)   \bline(0,0)(0,1)(150)
       \bline(0,150)(1,0)(150) \bline(150,0)(0,1)(150)
       \bezier{100}(0,0)(75,75)(150,150)
       \bezier{100}(85,0)(85,75)(85,150)
       \thicklines
       \bezier{200}(0,0)(40,35)(85,150)
       \bezier{200}(85,0)(115,45)(150,150)
       \put(-5,-8){$0$} \put(147,-8){$1$}
       \put(73,7){$c_\alpha$} %\put(100,-8){$X_2$}
       \put(-7,145){$1$}
       \bline(0,0)(1,0)(20) \bline(20,-3)(0,1)(3)
                            \put(7,-10){$\Delta_1$}
       \bline(80,0)(1,0)(20) \put(85,-10){$\Delta_{\jal}$}
       \bline(80,-3)(0,1)(3) \bline(100,-3)(0,1)(3)
      }}
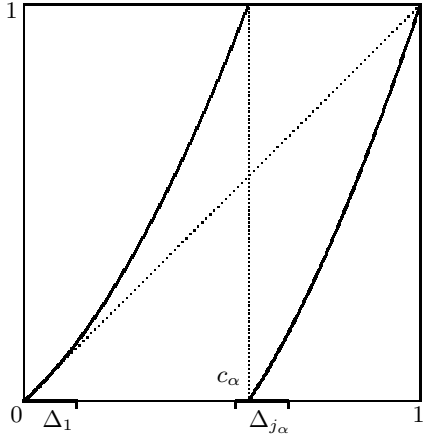
{The map $\map_\alpha$ and the positions of $c_\alpha$ and
$\Delta_1, \Delta_\jal$. \label{f:M-P}}
%%%%%%%%%%%%%%%%%%%%%%%%%%%%%%%%%%%%%%%%%%%%%%%%%%%%%%%%%%%%%%%%%

The idea of the proof of the items (a) - (c) is to make use
of the ``mass transfer'' between the intervals $\{\Delta_i\}$
under the action of the operator $\*{Q_\Delta} \*{\map_\alpha}$.
A sketch of the ``mass transfer'' together with the positions
of a few important intervals used in the proof are shown
in Fig.~\ref{fig-sketch}.

%%%%%%%%%%%%%%%%%%%%%%%%%%%%%%%%%%%%%%%%%%%%
%% Mass transfer
\Bfig(300,75)
      {\footnotesize{
       \put(10,0){\fbox{$\Delta_{1}$}}
       \put(55,0){\fbox{$\Delta_{2}$}}
       \put(80,0){{\bf\dots}}
       \put(100,0){\fbox{$\Delta_{\ell(2)}$}}
       \put(140,0){{\bf\dots}}
       \put(170,0){\fbox{$\Delta_{\ell(1)}$}}
       \put(230,0){{\bf\dots}}
       \put(260,0){\fbox{$\Delta_{\ell}$}}
       \put(290,0){{\bf\dots}}
       \bezier{200}(20,10)(125,65)(270,10)  %\ell --> 1
           \put(30,15){\vector(-2,-1){10}}  %----------
           \put(130,43){$p_{\ell,1}$}       %----------
       \bezier{20}(30,5)(38,5)(49,5)     %1 --> 2
           \put(49,5){\vector(1,0){6}} \put(35,10){$p_{1,2}$}
       \bezier{30}(130,5)(150,5)(170,5)  %\ell(2) --> \ell(1)
           \put(164,5){\vector(1,0){6}}
           \put(132,10){$p_{\ell(2),\ell(1)}$}
       \bezier{40}(200,5)(230,5)(260,5)  %\ell(1) --> \ell
           \put(254,5){\vector(1,0){6}}
           \put(220,10){$p_{\ell(1),\ell}$}
       \bezier{30}(15,10)(7,35)(0,3)       %\1 --> 1
       \bezier{30}(15,-5)(7,-30)( 0,3)
           \put(12,-14){\vector(1,3){3}} \put(-17,2){$p_{1,1}$}
      }}
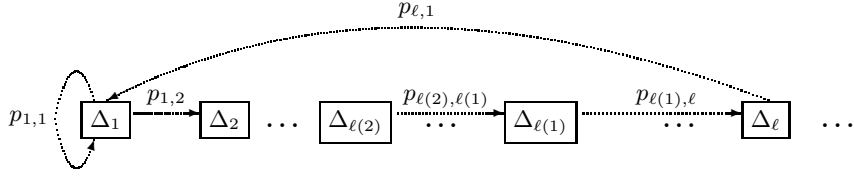
{Sketch of the connection between indices used in the proof.
\label{fig-sketch}}
%%%%%%%%%%%%%%%%%%%%%%%%%%%%%%%%%%%%%%%%%%%%%%%%%%%%%%%%%%%%%%%%%

Due to the monotonicity of the map $\mu_\Delta$ by
Lemma~\ref{l:key} for any $1\le i< j$ we have
$$ \frac{\pi_i}{|\Delta_i|} \ge \frac{\pi_j}{|\Delta_j|} $$
and hence
$$ \pi_i \ge \frac{|\Delta_i|}{|\Delta_j|}~\pi_j
         \ge \pi_j/K .$$
Recall that $|\Delta_i|/|\Delta_j|\le K$ for all $i,j$.

Denote by $z_1,z_2$ the unique solutions
to the equations %
\bea{z_1+z_1^{1+\alpha} \a= |\Delta_1| ,\\
     z_2+z_2^{1+\alpha} \a=
           {\rm min}\{|\Delta_1|+|\Delta_1|^{1+\alpha},
                          ~|\Delta_1|+|\Delta_2|\} .} %
Then $p_{1,1}=z_1/|\Delta_1|$ and
$p_{1,2}=(z_2-z_1)/|\Delta_1|$. Using that
$z_1^{1+\alpha}\le z_1\le|\Delta_1|<1$
for $\alpha\ge0$ we get $|\Delta_1|/2\le z_1\le|\Delta_1|$.

Observe now that for small enough $\delta>0$ one has
$|\Delta_1|^{1+\alpha}<|\Delta_2|$. Indeed,
$$ \frac{|\Delta_1|^{1+\alpha}}{|\Delta_2|}
 = \frac{|\Delta_1|}{|\Delta_2|}~|\Delta_1|^{\alpha}
\le K\delta^{\alpha} \toas{\delta\to0}0 .$$
Therefore
$$ {\rm min}\{|\Delta_1|+|\Delta_1|^{1+\alpha},
             ~|\Delta_1|+|\Delta_2|\}
 = |\Delta_1|+|\Delta_1|^{1+\alpha} $$
for $0<\delta\ll1$ and hence $z_2=|\Delta_1|$ and $p_{1,i}=0$ for
all $i>2$.

We have %
\beq{e:p12-up}{
   p_{1,2} = 1-p_{1,1} = 1-\frac{z_1}{|\Delta_1|}
             = \frac{z_1^{1+\alpha}}{|\Delta_1|}
             < \frac{|\Delta_1|^{1+\alpha}}{|\Delta_1|}
             = |\Delta_1|^{\alpha}
             \le \delta^{\alpha} .}%
To estimate $p_{1,2}$ from below we make use of that %
$\map_{\alpha}|\Delta_{1}|\le|\Delta_{1}| + |\Delta_{2}|$ if
$\delta\ll1$. %
\beq{e:p12-below}{
   p_{1,2} = \frac{z_2-z_1}{|\Delta_1|}
           = \frac{|\Delta_1|-z_1}{|\Delta_1|}
           = \frac{z_1^{1+\alpha}}{|\Delta_1|}
         \ge \frac{(|\Delta_1|/2)^{1+\alpha}}{|\Delta_1|}
         \ge 2^{-1-\alpha}\delta^{\alpha} .}%

Now we are ready to proceed with the proof of items (a) -- (c).

\bigskip

(a) Consider the interval $\Delta_\jal$ containing the right
endpoint $c_\alpha$ of the first interval of monotonicity of the
map $\map_\alpha$. Since the map $\map_\alpha$ is noncontracting
we have $p_{\jal+i,1}=0$ for all $i>K$. Indeed,
$\sum_{i=1}^{K}|\Delta_{\jal+i}|\ge|\Delta_{1}|$ and the map
$\map_{\alpha}$ is monotone on $\cup_{i=1}^{K}\Delta_{\jal+i}$.

By the definition of $\{\pi_{i}\}$
$$ \pi_{1} = p_{1,1}\pi_{1} +
             \sum_{i=0}^{K}p_{\jal+i,1}\pi_{\jal+i}$$
and hence %
\beq{e:p12p1-up}{
   p_{1,2} \pi_{1} = (1-p_{1,1}) \pi_{1}
                   = \sum_{i=0}^{K}p_{\jal+i,1}\pi_{\jal+i}
                 \le \sum_{i=0}^{K}\pi_{\jal+i}
                 \le (K+1)K\pi_{\jal} .}%

On the other hand, %
\beq{e:l-position}{
   \frac{c_\alpha}{\delta}
   \le \jal \le K\frac{c_\alpha}{\delta} .}%
Therefore
$$ 1 \ge \sum_{i=1}^{\jal}\pi_i
     \ge \frac{c_\alpha}{\delta} \cdot \frac{\pi_\jal}K
       = K^{-1}c_\alpha\delta^{-1}\pi_\jal ,$$
which implies %
\beaq{e:p1-up}{
  \pi_{1} \a\le \frac{(K+1)K}{p_{1,2}}~\pi_{\jal}
            \le (K+1)K2^{1+\alpha}\delta^{-\alpha}\cdot
                Kc_{\alpha}^{-1}\delta \cr
            \a= 2^{1+\alpha}(K+1)K^{2}c_{\alpha}^{-1}\delta^{1-\alpha} .}%

\smallskip

(b) Assume on the contrary that $\pi_1\le C_0\delta$ for some
$C_0<\infty$ and all $\delta\ll1$. Our aim is to show that this
assumption implies that $\mu_\Delta\toas{\delta\to0}\*{1_{\{0\}}}$
which is a contradiction. To demonstrate this convergence it is
enough to check that for any $z\in(0,1]$ we have
$\mu_\Delta([z,1])\toas{\delta\to0}0$.

Denote %
$$ \ell := \function{\jal   &\mbox{if } p_{\jal,1}\ge p_{\jal+1,1} \\
                     \jal+1 &\mbox{otherwise}}$$ %
and %
\beaq{e:lipt}{
  \lal \a:= \sup_{k\ge0}
            \map_{\alpha}'((\map_{\alpha}|[0,c_{\alpha}])^{-k}
                                    c_{\alpha}+\delta) \cr
         \a=  \max_{k\in\{0,1\}}
            \map_{\alpha}'((\map_{\alpha}|[0,c_{\alpha}])^{-k}
                                   c_{\alpha}+\delta) } %
due to the convexity of the first branch of the map $\map_\alpha$.
Here $(\map_{\alpha}|A)^{-1}x:=\map_{\alpha}^{-1}x\cap A$ and
$\delta$ is assumed to be small enough.

By the construction the value $p_{i,j}|\Delta_{i}|$ is equal to the
Lebesgue measure of the part of $\Delta_{i}$ which is mapped into
$\Delta_{j}$ by $\map_{\alpha}$. Therefore
$$ \lal(p_{\jal,1}|\Delta_{\jal}| + p_{\jal+1,1}|\Delta_{\jal+1}|)
   \ge |\Delta_{1}| .$$
Hence %
$$ \max_{i\in\{0,1\}}\{p_{\jal+i,1}\}\cdot
   \max_{i\in\{0,1\}}\{|\Delta_{\jal+i}|\}
   \ge\max_{i\in\{0,1\}}\{p_{\jal+i,1}|\Delta_{\jal+i}|\}
   \ge \frac{|\Delta_{1}|}{2\lal} $$
and thus %
\beq{e:pl1}{
 p_{\ell,1}\ge(2K\lal)^{-1} .}%

Considering the ``mass transfer'' between intervals $\Delta_1$ and
$\Delta_\ell$ and using the estimate for $p_{1,2}$ from above
(\ref{e:p12-up}) we get %
\beq{e:pil}{
   \pi_\ell \le \frac{p_{1,2}}{p_{\ell,1}}~\pi_1
            \le \delta^{\alpha}~2K\lal~\pi_1
            \le 2K\lal~C_0\delta^{1+\alpha} }%
by the assumption on $\pi_1$.

Similarly to the definition of the index $\ell$, one defines by
induction a sequence of indices $\{\ell(t)\}_{t\ge0}$ as follows.
We set $\ell(0):=\ell$ and
$\hat\ell(t):=\min\{j:~p_{j,\ell(t-1)}>0\}$ for $t\ge1$. Then for
$t\ge1$ we set
$$ \ell(t) := \function{\hat\ell(t)   &\mbox{if }
                p_{\hat\ell(t),\ell(t-1)}\ge p_{\hat\ell(t)+1,\ell(t-1)}\\
                      \hat\ell(t) +1 &\mbox{otherwise} .}$$ %
Using the same argument as above one estimates the transition
probabilities from below as follows:
$$ p_{\ell(t),\ell(t-1)}\ge(2K\lal)^{-1} .$$
By the construction %
$\map_\alpha\Delta_{\ell(t)} \cap \Delta_{\ell(t)}=\emptyset$ %
for small enough $\delta>0$ and any $t\ge0$. Therefore
$$ \pi_{\ell(t-1)} \ge p_{\ell(t),\ell(t-1)}\pi_{\ell(t)} .$$
for any $t\ge1$. Hence by (\ref{e:pil})
$$ \pi_{\ell(t)} \le
   \left(\prod_{i=1}^{t}p_{\ell(i),\ell(i-1)}\right)^{-1}\pi_{\ell}
   \le (2K\lal)^{t+1}C_{0}\delta^{1+\alpha} .$$

Consider a sequence of points $\{\beta_t\}_{t\ge0}$ from the interval
$(0,c_\alpha]$ such that
$$ \beta_0:=c_\alpha, ~~ \map_\alpha\beta_t=\beta_{t-1}
                   ~~{\rm for}~~ t\ge1 .$$
This sequence converges to zero as $n\to\infty$ and for this
specific map one even can get an asymptotic formula
$\beta_n\approx Cn^{-1/\alpha}$ for $n\to\infty$ (see \cite{LVS}).
For our aim it is enough to observe that for any
$z\in(0,c_\alpha]$ there exists a finite index $t_z$ such that
$z\ge\beta_{t_z}$. On the other hand,
$$ x_{t+1} \le \beta_t \le x_t + 2\delta $$
for any $t\ge0$ and any pair of points $x_t\in\Delta_{\ell(t)},
~~x_{t+1}\in\Delta_{\ell(t+1)}$, provided $\delta>0$ is small
enough. Therefore
$$ z\ge\beta_{t_z} \ge x_{t_z+1} $$
for any $x_{t_z+1}\in\Delta_{\ell(t_z+1)}$.

Making use of the assumption $\pi_1\le C_0\delta$ we obtain %
\bea{ \mu_\Delta([z,1]) \a\le \mu_\Delta([\beta_{t_z},1])
       \le \sum_{i\ge \ell(t_z+1)} \pi_i \\
       \a\le \frac{1}{\delta/K}~K\pi_{\ell(t_z+1)}
       \le \frac{K^2}{\delta}~(2K\lal)^{t_{z}+2}
           C_{0}\delta^{1+\alpha} \\
       \a= K^2(2K\lal)^{t_{z}+2}C_0\delta^{\alpha}
       \toas{\delta\to0}0 }%
for any $\alpha>0$. This proves (b).

Observe now that item (b) together with item (a) and the
monotonicity of the measures under study implies that the measure
$\mu_\Delta$ does not converge to the Dirac measure at the
origin. Note however that this is not enough to prove the
convergence to the absolute continuous SRB measure of
$\map_\alpha$ existing for $0<\alpha<1$.

\smallskip

(c) The proof of the remaining part is very similar to the
previous one except that we do not need to make any additional
assumptions. Using the notation introduced in the
proof of item (b) we get %
\bea{ \mu_\Delta([z,1]) \a\le \mu_\Delta([\beta_{t_z},1])
       \le \sum_{i\ge \ell(t_z+1)} \pi_i \\
       \a\le \frac{1}{\delta/K}~K\pi_{\ell(t_z+1)}
       \le \frac{K^2}{\delta}~(2K\lal)^{t_{z}+2}
           \frac{p_{1,2}}{p_{\ell,1}}\pi_1 \\
       \a\le \frac{K^2}{\delta}~(2K\lal)^{t_{z}+2}
             \delta^{\alpha}\pi_1 \\
       \a\le K^2(2K\lal)^{t_{z}+2}\delta^{\alpha-1}
       \toas{\delta\to0}0 .}%
because $\pi_1\le1$ and $\alpha>1$. Now since $z\in(0,c_\alpha]$
is arbitrary this implies that $\mu_\Delta\toas{\delta\to0}
\*{1_{\{0\}}}$. \qed

\section{Generalizations} \label{s:generalization}

A close look to the proof of Theorem~\ref{t:main} in the previous
Section shows that we were using very few specific properties of
the Mannevile-Pomeau map $\map_\alpha$ and while the fact
that this map belongs to the family of piecewise convex maps
${\cal T}$ introduced in Section~\ref{s:monotonic} is used heavily.
The aim of this Section is to demonstrate that indeed adding a few
assumptions to the definition of the family ${\cal T}$ one can
prove the result of the same sort as Theorem~\ref{t:main}.

\begin{theorem}\label{t:generalization}
Let $\map\in{\cal T}$ and let it satisfy the following
assumptions %
\begin{enumerate}
\item[(i)] $\map x = x + Cx^{1+\alpha} + o(x^{1+\alpha})$
      as $x\to0$ with $\alpha>0$,
\item[(ii)] $|\map x  -  \map y|\ge|x-y|$ for all $x,y\in X$
      such that $|x-y|\ll1$,
\item[(iii)] {\rm Card}$\{\map^{-1}x\}\le M<\infty$ for any $x\in X$. %
\end{enumerate}
Then all the claims made in Theorem~\ref{t:main} remain valid in
this setting.
\end{theorem}

Observe that the map $\map$ needs not to be Markov and the number
of branches of the inverse map $\map^{-1}$ is arbitrary (but
finite).

\proof The scheme of the proof is exactly the same as in the case
of Theorem~\ref{t:main} and we explain only how to overcome
difficulties related to our more general setup.

The main difference between the situations considered in these
two Theorems is that the source of the ``mass transfer'' to the
interval $\Delta_1$ is no longer restricted to the beginning of
the interval $X_2$ of the corresponding special partition, namely
to the interval $\Delta_\jal$ defined in the proof of
Theorem~\ref{t:main}. Moreover, we need to make the
assumption~(ii) that the map is noncontracting (since in general
a map from the family ${\cal T}$ needs not to satisfy this
property).

In the present setting in the beginning of each element of the
special partition there is an interval of the partition $\Delta$
playing the same role as $\Delta_\jal$. Nevertheless the ``mass
transfer'' from these additional sources may only enlarge the
amount arriving to $\Delta_1$ and $\Delta_{\ell(t)}$ and thus do
not change the estimates which we use in the proof of the items
(b) and (c).

To take care about these additional sources in the proof of the
item (a) we make use of the assumption (iii) which enables us to
estimate the number of these sources and to write a variant of
the inequality (\ref{e:p12p1-up}) as follows %
\bea{p_{1,2} \pi_{1} \a= (1-p_{1,1}) \pi_{1} \cr
                 \a\le KM\sum_{i=0}^{K}\pi_{\jal+i}
                   \le (K+1)K^{2}M\pi_{\jal} .}%
Observe that here we use heavily the monotonicity of the invariant
distribution.

Applying this estimate instead of (\ref{e:p12p1-up}) and following
the same arguments as in the proof of the item (a) of
Theorem~\ref{t:main} one gets %
$$ \pi_{1} \le
   2^{1+\alpha}(K+1)K^{3}Mc_{\alpha}^{-1}\delta^{1-\alpha} .$$

It remains to discuss the estimates of $p_{1,2}$ which for
$\delta$ small enough depend only on the behavior of the map
$\map$ in a small neighborhood of the origin. Using the same
notation as in the proof of Theorem~\ref{t:main} consider the
unique solutions to the equations %
\bea{z_1+Cz_1^{1+\alpha} \a= |\Delta_1| ,\\
     z_2+Cz_2^{1+\alpha} \a=
           {\rm min}\{|\Delta_1|+C|\Delta_1|^{1+\alpha},
                          ~|\Delta_1|+|\Delta_2|\} .} %
Then for small enough $0<\delta\ll1$ we get
$$ p_{1,2} \le 1-p_{1,1} = 1-\frac{z_1}{|\Delta_1|}
             = \frac{Cz_1^{1+\alpha}}{|\Delta_1|}
             < C\frac{|\Delta_1|^{1+\alpha}}{|\Delta_1|}
             \le C\delta^{\alpha} $$
and
$$ p_{1,2} = \frac{z_2-z_1}{|\Delta_1|}
           = \frac{Cz_1^{1+\alpha}}{|\Delta_1|}
         \ge C\frac{(|\Delta_1|/(2C))^{1+\alpha}}{|\Delta_1|}
         \ge 2^{-1-\alpha}C^{-\alpha}\delta^{\alpha} .$$
Taking into account that the term $o(x^{1+\alpha})$ in (i) gives
only higher order corrections to the estimates above one applies
directly all further arguments used in the proof of
Theorem~\ref{t:main} in the present setting as well. \qed

\section*{Acknowledgments}
The author is grateful to anonymous referees for helpful comments
and suggestions which improved the quality of the paper.

\small%\newpage
%%%%%%%%%%%

\end{document}